\documentclass[a4paper,11pt]{article}

\usepackage{natbib}
\usepackage{amssymb,amsmath,amsfonts,stmaryrd}
\usepackage[latin1]{inputenc}
\usepackage{color,rotating}
\usepackage{authblk}

\newtheorem{theorem}{Theorem}

\newcommand{\esp}{\mathrm{E}}
\newcommand{\matr}[1]{\boldsymbol{#1}}
\newcommand{\tr}{\mathrm{tr}}
\newcommand{\transp}{^{\sf t}}
\newcommand{\ud}{\mathrm{d}}
\newcommand{\var}{\mathrm{Var}}
\newcommand{\vect}[1]{\boldsymbol{#1}}

\newcommand{\idep}{i_d}
\newcommand{\jdep}{j_d}

\newcommand{\gcite}{\citep}
\newcommand{\gcitet}{\citet}
\newcommand{\gcitep}{\citealp}

\title{Cumulants of multiinformation density in the case of a multivariate normal distribution}
\author[1,2]{Guillaume Marrelec}
\author[1,2]{Alain Giron}

\affil[1]{Laboratoire d'imagerie biom{\'e}dicale (LIB), Sorbonne Universit{\'e}, CNRS, INSERM, F-75006, Paris, France. email: firstname.lastname@inserm.fr.}
\affil[2]{Centre de recherches et d'études en sciences des interactions (CR{\'E}SI), Center for Interaction Science (CIS), F-75006, Paris, France}

\begin{document}

\maketitle

\begin{abstract}
We consider a generalization of information density to a partitioning into $N \geq 2$ subvectors. We calculate its cumulant-generating function and its cumulants, showing that these quantities are only a function of all the regression coefficients associated with the partitioning.
\par
\
\par
\noindent
\textit{Keywords:} dependence; information; cumulant-generating function; cumulants; mutual information; multiinformation
\end{abstract}

\section{Introduction}

%\comment{Attention: ``Articles will be limited to six journal pages (about 13 double-space typed pages) including references and figures''}

Let $\vect{X}$ be a multivariate normal random variable with distribution $f ( \vect{x} )$, and $( \vect{X}_1, \vect{X}_2 )$ a partitioning of $\vect{X}$ into 2 subvectors with corresponding marginals $f_1 ( \vect{x}_1 )$ and $f_2 (\vect{x}_2 )$. The information density relative to $\vect{X}$ and the partitioning $( \vect{X}_1, \vect{X}_2 )$ is the random variable defined as \gcite[\S17.1]{Polyanskiy-2017}
\begin{equation} \label{eq:i:def}
 \idep ( \vect{X}; \vect{X}_1, \vect{X}_2 ) = \ln \frac{ f ( \vect{X}_1, \vect{X}_2 ) } { f_1 ( \vect{X}_1 ) \, f_2 ( \vect{X}_2 ) }.
\end{equation}
One of the key features of information density is that its expectation yields mutual information (\gcitep[Chap.~1, \S2]{Kullback-1968}; \gcitep[\S17.1]{Polyanskiy-2017}). In the present paper, we consider a partitioning of $\vect{X}$ into $N \geq 2$ subvectors $( \vect{X}_1, \dots, \vect{X}_N )$ with corresponding marginals $f_1, \dots, f_N$ and define \emph{multiinformation density} as
$$\idep ( \vect{X}; \vect{X}_1, \dots, \vect{X}_N ) = \ln \left[ \frac{ f ( \vect{X} ) } { \prod_{ n = 1 } ^ N f_n ( \vect{X}_n ) } \right].$$
The mean of this quantity 
$$I ( \vect{X}; \vect{X}_1, \dots, \vect{X}_N ) = \int f ( \vect{x} ) \, \ln \left[ \frac{ f ( \vect{x} ) } { \prod_{ n = 1 } ^ N f_n ( \vect{x}_n ) } \right] \, \ud \vect{x}$$
is itself a generalization of mutual information known under different names: total correlation \gcite{Watanabe-1960}, multivariate constraint \gcite{Garner-1962}, $\delta$ \gcite{Joe-1989b}, or multiinformation \gcite{Studeny-1998b}. Our interest in $i_d$ is driven by its close connection with mutual independence. Indeed, when the $\vect{X}_n$'s are mutually independent, multiinformation is classically equal to 0, but we also have $i_d \equiv 0$ and $\var ( i_d ) = 0$ (see Appendix~\ref{ann:id}), yielding other statistical markers of independence. By contrast, dependence between the $\vect{X}_n$'s is a multivariate phenomenon that multiinformation, as a one-dimensional measure, can only partially quantify. We expect $i_d$ to give a more detailed characterization of dependence, e.g., through its moments or cumulants.
\par
We here focus on multivariate normal distributions. The family of multivariate normal distribution with given mean $\vect{\mu}$ can be parameterized by either a covariance matrix or a concentration/precision (i.e., inverse covariance) matrix. Either parameterization shows multivariate distributions according to a certain perspective and emphasizes different features (e.g., Markov properties for the concentration matrix). In this  context, we wished to investigate the existence of a natural way of parameterizing dependencies, i.e., a parameter that would emphasize the dependence properties of the distribution.
\par
The core of the present paper is the following theorem:
\begin{theorem}
 Let $\vect{X}$ be a $d$-dimensional variable following a multivariate normal distribution with mean $\vect{\mu}$ and covariance matrix $\matr{\Sigma}$.  Partition $\vect{X}$ into $N$ subvectors $( \vect{X}_1, \dots, \vect{X}_N )$, and set $\idep$ the corresponding multiinformation density. Then the cumulant-generating function of $i_d$ is given by
 \begin{equation} \label{eq:def:cgf}
 \ln \esp \left( e ^ { t \idep } \right) = t \, I ( \vect{X}_1; \dots; \vect{X}_N ) - \frac{1}{2} \ln | \matr{I}_d - t \matr{\Gamma} |,
\end{equation}
 where
 $$\matr{\Gamma} = \matr{\Sigma} \, \mathrm{diag} ( \matr{\Sigma}_{ 1 1 }, \dots, \matr{\Sigma}_{ N N } ) ^ { - 1 } - \matr{I}_d$$
 is the block matrix whose diagonal blocks are equal to $\matr{0}$ and where each off-diagonal block $( m, n )$ is the matrix of regression coefficients of $\vect{X}_m$ on $\vect{x}_n$
 \begin{equation} \label{eq:def:regcoeff}
 \matr{\Gamma}_{ m | n } = \matr{\Sigma}_{ m n } \matr{\Sigma}_{ n n } ^ { - 1 }, \qquad m \neq n.
\end{equation}
 The cumulants of $i_d$ are given by
 $$\kappa_1 ( \idep ) = I ( \vect{X}_1; \dots, \vect{X}_N ) \qquad \mbox{and} \qquad \kappa_l ( \idep ) = \frac{ ( l - 1 ) ! } { 2 } \, \tr \left( \matr{\Gamma} ^ l \right), \ l \geq 2.$$
\end{theorem}
This theorem is proved in Section~\ref{s:preuve}. In Section~\ref{s:cons}, we investigate some consequences of this result. Section~\ref{s:disc} is devoted to the discussion.

\section{Proof of theorem} \label{s:preuve}

\subsection{Cumulant-generating function}

We partition $\vect{\mu}$ and $\matr{\Sigma}$ in accordance with the partitioning of $\vect{X}$, so that $\vect{\mu}_n$ is the expectation of $\vect{X}_n$ and $\matr{\Sigma}_{ m n }$ the matrix of covariances between $\vect{X}_m$ and $\vect{X}_n$. Multiinformation between the $\vect{X}_n$'s yields
$$I ( \vect{X}_1; \dots; \vect{X}_N ) = \frac{ 1 } { 2 } \ln \frac{ \prod_{ n = 1 } ^ N | \matr{\Sigma}_{ n n } | } { | \matr{\Sigma} | }.$$
From there, we can express multiinformation density as
\begin{equation} \label{eq:def:idep-jdep}
 \idep = I ( \vect{X}_1; \dots; \vect{X}_N ) + \jdep,
\end{equation}
where $\jdep$ is defined as
\begin{eqnarray} \label{eq:def:jdep}
 \jdep %& = & \sum_{ n = 1 } ^ N ( \vect{X}_n - \vect{\mu}_n ) \transp \matr{\Sigma}_{ n n } ^ { - 1 }( \vect{X}_n - \vect{\mu}_n ) - \frac{ 1 } { 2 } ( \vect{X} - \vect{\mu} ) \transp \matr{\Sigma} ^ { - 1 } ( \vect{X} - \vect{\mu} ) \\
 & = & \frac{ 1 } { 2 } ( \vect{X} - \vect{\mu} ) \transp \matr{\Phi} ( \vect{X} - \vect{\mu} ).
\end{eqnarray}
and $\matr{\Phi}$ as
$$\matr{\Phi} = \mathrm{diag} ( \matr{\Sigma}_{ 1 1 }, \dots, \matr{\Sigma}_{ N N } ) ^ { - 1 } - \matr{\Sigma} ^ { - 1 }.$$
Here, $\mathrm{diag} ( \matr{\Sigma}_{ 1 1 }, \dots, \matr{\Sigma}_{ N N } )$ stands for the block-diagonal matrix with diagonal blocks equal to the $\matr{\Sigma}_{ n n }$'s. The moment-generating function of $\jdep$ yields
\begin{eqnarray*}
 \esp \left( e ^ { t \jdep } \right) % & = & \int e ^ { \frac{ t } { 2 } ( \vect{x} - \vect{\mu} ) \transp \matr{\Phi} ( \vect{x} - \vect{\mu} ) } \, f ( \vect{x} ) \, \ud x \\
 & = & \int ( 2 \pi ) ^ { - \frac{ d } { 2 } } | \matr{\Sigma} | ^ { - \frac{ 1 } { 2 } } e ^ { - \frac{ 1 } { 2 } ( \vect{x} - \vect{\mu} ) \transp \left( \matr{\Sigma} ^ { - 1 } - t \matr{\Phi} \right) ( \vect{x} - \vect{\mu} ) } \, \ud \vect{x}.
\end{eqnarray*}
Since it can be shown that $\matr{\Sigma} ^ { - 1 } - t \matr{\Phi}$ is positive definite at least in a neighborhood of $t = 0$ (see Appendix~\ref{an:pd}), the integrand is proportional to a multivariate normal distribution with mean $\vect{\mu}$ and covariance matrix $\matr{\Sigma} ^ { - 1 } - t \matr{\Phi}$. Integration with respect to $\vect{x}$ therefore yields
$$\esp \left( e ^ { t \jdep } \right) = \frac{ | \matr{\Sigma} ^ { - 1 } | ^ { \frac{ 1 } { 2 } } } { | \matr{\Sigma} ^ { - 1 } - t \matr{\Phi} | ^ { \frac{ 1 } { 2 } } } = | \matr{I}_d - t \matr{\Gamma} | ^ { - \frac{ 1 } { 2 } }$$
and
\begin{equation} \label{eq:def:cgf}
 \ln \esp \left( e ^ { t \jdep } \right) = - \frac{1}{2} \ln | \matr{I}_d - t \matr{\Gamma} |,
\end{equation}
% \begin{eqnarray*}
%  \esp \left( e ^ { t \jdep } \right) & = & \frac{ | \matr{\Sigma} ^ { - 1 } | ^ { \frac{ 1 } { 2 } } } { | \matr{\Sigma} ^ { - 1 } - t \matr{\Phi} | ^ { \frac{ 1 } { 2 } } } \\
%  & = & | \matr{I}_d - t \matr{\Gamma} | ^ { - \frac{ 1 } { 2 } }, 
% \end{eqnarray*}
where $\matr{I}_d$ is the $d$-by-$d$ unit matrix and $\matr{\Gamma} = \matr{\Sigma} \matr{\Phi}$ the block matrix whose diagonal blocks are equal to $\matr{0}$ and where each nondiagonal block $( m, n )$ is the matrix of regression coefficients of $\vect{X}_m$ on $\vect{x}_n$ given by \gcite[Definition~2.5.1]{Anderson_TW-2003}
\begin{equation} \label{eq:def:regcoeff}
 \matr{\Gamma}_{ m | n } = \matr{\Sigma}_{ m n } \matr{\Sigma}_{ n n } ^ { - 1 }, \qquad m \neq n.
\end{equation}

\subsection{Cumulants}

The cumulants of $\idep$ can be calculated in closed form from those of $\jdep$ and Equation~\eqref{eq:def:idep-jdep} by noting that the first cumulant, $\kappa_1 ( \idep )$, is shift-equivariant, while the others, $\kappa_i ( \idep )$ for $i \geq 2$, are shift invariant \gcite[\S3.13]{Kendall-1945}. This leads to
\begin{equation} \label{eq:I:cum:J}
 \left\{ \begin{array}{ccll}
 \kappa_1 ( \idep ) & = & I ( \vect{X}_1; \dots; \vect{X}_N ) + \kappa_1 ( \jdep ) & \\
 \kappa_l ( \idep ) & = & \kappa_l ( \jdep ), & l \geq 2.
\end{array} \right.
\end{equation}
Now, the cumulants of $\jdep$ can be easily computed as follows. Using the fact that $| \matr{A} | = e ^ { \tr [ \ln ( \matr{A} ) ] }$ \gcite{Higham-2007}, which, for a positive definite matrix, can be expressed as $\ln | \matr{A} | = \tr [ \ln ( \matr{A} ) ]$, we have from Equation~\eqref{eq:def:cgf}
$$\ln \esp \left( e ^ { t \jdep } \right) = - \frac{ 1 } { 2 } \tr \left[ \ln ( \matr{I}_d - t \matr{\Gamma} ) \right].$$
For $t$ sufficiently small, we can perform a Taylor expansion of the log function around $\matr{I}_d$ \gcite[Eq.~4.1.24]{Abramowitz-1972}, leading to
\begin{eqnarray*}
 \ln \esp \left( e ^ { t \jdep } \right) & = & \frac{ 1 } { 2 } \tr \left[ \sum_{ il = 1 } ^ { \infty } \frac{ ( t \matr{\Gamma} ) ^ l  } { l } \right] \\
 & = & \sum_{ l = 1 } ^ { \infty } \frac{ t ^ l } { 2 l } \, \tr ( \matr{\Gamma} ^ l ).
\end{eqnarray*}
Identification with the decomposition of the same function in terms of cumulants \gcite[\S3.12]{Kendall-1945}
$$\ln \esp \left( e ^ { t \jdep } \right) = \sum_{ l = 1 } ^ { \infty } \kappa_l \frac{ t ^ l } { l ! },$$
yields for the cumulants of $\jdep$
\begin{equation} \label{eq:J:cum}
 \kappa_l ( \jdep ) = \frac{ ( l - 1 ) ! } { 2 } \, \tr ( \matr{\Gamma} ^ l ).
\end{equation}
The same result could have been reached by using the fact that $\jdep$ is a quadratic function of a multidimensional normal variate $\vect{x}$, as evidenced in Equation~\eqref{eq:def:jdep}, together with the expression of the cumulants of such functions \gcite[Lemma~2]{Magnus-1986}.
\par
The cumulants of $\idep$ therefore yield
$$\kappa_1 ( \idep ) = I ( \vect{X}_1; \dots, \vect{X}_N ),$$
as expected, since the first cumulant is also the mean \gcite[\S3.14]{Kendall-1945}, and, for $l \geq 2$,
$$\kappa_l ( \idep ) = \frac{ ( l - 1 ) ! } { 2 } \, \tr \left( \matr{\Gamma} ^ l \right).$$
In particular, the variance, which is equal to the second cumulant \gcite[\S3.14]{Kendall-1945} is given by
\begin{equation} \label{eq:idep:var}
 \var ( \idep ) = \kappa_2 ( \idep ) = \sum_{ 1 \leq m < n \leq N } \tr \left( \matr{\Sigma}_{ m n } \matr{\Sigma}_{ m m } ^ { - 1 } \matr{\Sigma}_{ m n } \matr{\Sigma}_{ n n } ^ { - 1 } \right).
\end{equation}
In the even more particular case where all subvectors are 1-dimensional, we have
$$\var ( \idep ) =  \sum_{ 1 \leq m < n \leq N } \rho_{ m n } ^ 2,$$
where $\rho_{mn}$ is the usual correlation coefficient between $X_m$ and $X_n$.

\section{Consequences} \label{s:cons}

We here investigate some consequences of the previous results: the particular case of partitioning into two subvectors, the irrelevance of the variances, and a graphical interpretation.

\subsection{Partitioning into two subvectors}

In the particular case where $N = 2$, multiinformation boils down to mutual information. The various powers of $\matr{\Gamma}$ can easily be calculated, yielding
$$\matr{\Gamma} ^ l = \begin{pmatrix}
 \matr{0} & \matr{\chi}_l \\
 \matr{\Upsilon}_l & 0
\end{pmatrix} \qquad \mbox{for $l$ odd}$$
and 
$$\matr{\Gamma} ^ l = \begin{pmatrix}
 \matr{\chi}_l & \matr{0} \\
 \matr{0} & \matr{\Upsilon}_l
\end{pmatrix} \qquad \mbox{for $l$ even},$$
where we set
$$\left\{ \begin{array}{ccl}
 \matr{\chi}_l & = & \underbrace{\matr{\Gamma}_{ 1 | 2 } \matr{\Gamma}_{ 2 | 1 } \dots }_{\mbox{$l$ factors}} \\
 \matr{\Upsilon}_l & = & \overbrace{\matr{\Gamma}_{ 2 | 1 } \matr{\Gamma}_{ 1 | 2 } \dots },
\end{array} \right.$$
with the relationship that $\matr{\chi}_l = \matr{\Gamma}_{ 1 | 2 } \matr{\Upsilon}_{ l - 1 }$ and $\matr{\Upsilon}_l = \matr{\Gamma}_{ 2 | 1 } \matr{\chi}_{ l - 1 }$. For $l$ odd, the trace of $\matr{\Gamma} ^ l$ is equal to 0; for $l$ even, it is equal to $\tr ( \matr{\chi}_l ) + \tr ( \matr{\Upsilon}_l )$, which can alternatively be expressed as
$$2 \, \tr ( \matr{\chi}_l ) = 2 \, \tr \left[ ( \matr{\Gamma}_{ 1 | 2 } \, \matr{\Gamma}_{ 2 | 1 } ) ^ \frac{ l } { 2 } \right]$$
or
$$2 \, \tr ( \matr{\Upsilon}_l ) = 2 \, \tr \left[ ( \matr{\Gamma}_{ 2 | 1 } \, \matr{\Gamma}_{ 1 | 2 } ) ^ \frac{ l } { 2 } \right].$$
In particular, the variance of $\idep$ is equal to
\begin{equation} \label{eq:2:var}
 \var ( \idep ) = \kappa_2 ( \idep ) = \tr \left( \matr{\Sigma}_{ 12 } \matr{\Sigma}_{ 22 } ^ { - 1 } \matr{\Sigma}_{ 21 } \matr{\Sigma}_{ 11 } ^ { - 1 } \right).
\end{equation}
This quantity, which is a particular case of Equation~\eqref{eq:idep:var}, was introduced by \gcitet{Jupp-1980} as an extension of the classical correlation coefficient in the case of multidimensional variates, with application to directed data \gcite[\S11.2]{Mardia-2000}. It is the sum of the squared canonical correlation coefficients between $\vect{X}_1$ and $\vect{X}_2$ (\gcitep[Chap.~12]{Anderson_TW-2003}; \gcitep{Jupp-1980}).
\par
If we furthermore assume that $X_1$ is a 1-dimensional vector, the cumulant-generating function yields
$$\ln \esp \left( e ^ { t \idep } \right) = - \frac{ t } { 2 } \ln ( 1 - \overline{R} ^ 2 ) - \ln ( 1 - t ^ 2 \overline{R} ^ 2 ),$$
where $\overline{R}$ is the multiple correlation coefficient \gcite[\S~2.5.2]{Anderson_TW-2003}
$$\overline{R} ^ 2 = \frac{ \vect{\Sigma}_{ 12 } \matr{\Sigma}_{ 22 } ^ { - 1 } \vect{\Sigma}_{ 12 } \transp } { \sigma_1 ^ 2 }.$$
The cumulants are equal to 0 for $l$ odd and to $\kappa_l ( \idep ) = ( l - 1 )! \, ( \overline{R} ^ 2 ) ^ { \frac{ l } { 2 } }$ for $l$ even. In particular, the variance reads $\var ( \idep ) = \overline{R} ^ 2$.
\par
Finally, if both $X_1$ and $X_2$ are assumed to be 1-dimensional vectors with correlation coefficient $\rho$, the cumulant-generating function reads
$$\ln \esp \left( e ^ { t \idep } \right) = - \frac{ t } { 2 } \ln ( 1 - \rho ^ 2 ) - \ln ( 1 - t ^ 2 \rho ^ 2 ),$$
with cumulants equal to 0 for $l$ odd and $\kappa_l ( \idep ) = ( l - 1 )! \, | \rho | ^ l$ for $l$ even. In particular, the variance yields $\var ( \idep ) = \rho ^ 2$.

\subsection{Irrelevance of variances}

Mutual information and multiinformation are both quantities that do not depend on the variance coefficients $\Sigma_{kk}$'s. This result can be generalized to the cumulant-generating function of $\idep$, and hence, its distribution, moments and cumulants. Indeed, let
$$\matr{\Sigma} = \matr{\Delta} \matr{R} \matr{\Delta}$$
be the decomposition of the covariance matrix $\matr{\Sigma}$ such that $\matr{\Delta} = ( \Delta_{ k l } )$ is a diagonal matrix with $\Delta_{ k k } = \sqrt{ \Sigma_{ k k } }$ and $\matr{R} = ( R_{ k l } )$ is the correlation matrix with $R_{ k l } = \Sigma_{ k l } / \sqrt{ \Sigma_{ k k } \Sigma_{ l l } }$. Block multiplication of $\matr{\Delta} \matr{R} \matr{\Delta}$ shows that we have
$\matr{\Sigma}_{ m n } = \matr{\Delta}_{ m m } \matr{R}_{ m n } \matr{\Delta}_{ n n }$ for any block $\matr{\Sigma}_{ m n }$ of $\matr{\Sigma}$. Using the fact that $( \matr{A} \matr{B} ) ^ { - 1 } = \matr{B} ^ { - 1 } \matr{A} ^ { - 1 }$ for any two invertible square matrices $\matr{A}$ and $\matr{B}$, we have
$$\matr{\Gamma}_{ m | n } = \matr{\Delta}_{ m m } \matr{R}_{ m n } \matr{R}_{ n n } ^ { - 1 } \matr{\Delta}_{ n n } ^ { - 1 } = \matr{\Delta}_{ m m } \tilde{\matr{\Gamma}}_{ m | n } \matr{\Delta}_{ n n } ^ { - 1 },$$
where $\tilde{\matr{\Gamma}}_{ m | n }$ is the matrix of regression coefficients obtained by application of Equation~\eqref{eq:def:regcoeff} to the correlation matrix $\matr{R}$ instead of the covariance matrix $\matr{\Sigma}$. This result shows that $\matr{\Gamma} = ( \matr{\Gamma}_{ m | n } )$ can be factorized into
$$\matr{\Gamma} = \matr{\Delta} \tilde{\matr{\Gamma}} \matr{\Delta} ^ { - 1 },$$
where we set $\tilde{\matr{\Gamma}} = ( \tilde{\matr{\Gamma}}_{ m | n } )$. This takes us to
\begin{eqnarray*}
 | \matr{I}_d - t \matr{\Gamma} | & = & | \matr{\Delta} ( I_d - t \tilde{\matr{\Gamma}} ) \matr{\Delta} ^ { - 1 } | \\
 & = & | I_d - t \tilde{\matr{\Gamma}} |.
\end{eqnarray*}
As a conclusion, we have that the cumulant-generating function of a multivariate distribution with covariance matrix $\matr{\Sigma}$ is the same as the cumulant-generating function of a multivariate distribution with covariance matrix $\matr{R}$, where $\matr{R}$ is the correlation matrix associated with $\matr{\Sigma}$. This is a translation of the fact that $\idep$ does not depend on the variance coefficients.

\subsection{Graphical interpretation} \label{ss:gi}

\begin{figure}[!htbp]
 \centering
 \includegraphics[width=5cm]{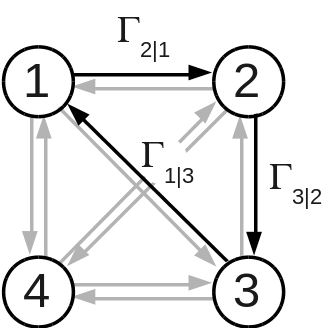}
 \caption{{\bf Graphical interpretation of $\tr ( \matr{\Gamma} ^ l )$.} We consider the case $N = 4$ and $l = 3$. The directed 3-loop $lp = (1 \to 2 \to 3 \to 1 )$ is represented with dark arrows. The value of $\tau$ on this loop is equal to $\tau ( l ) = \tr ( \matr{\Gamma}_{ 1 | 3 } \matr{\Gamma}_{ 3 | 2 } \matr{\Gamma}_{ 2 | 1 } )$. $\tr ( \matr{\Gamma} ^ 3 )$ is obtained by summing $\tau ( p )$ over all 3-loops.}
 \label{fig}
\end{figure}

While simplification of the expression for the cumulants through the explicit calculation of $\matr{\Gamma} ^ l$ is challenging in the general case, one can resort to a graphical interpretation of this matrix. Note first that the block $( n, m )$ of $\matr{\Gamma} ^ l$ is given by
$$( \matr{\Gamma} ^ l )_{ m n } = \sum_{ q_1, \dots, q_{ l - 1 } } \underbrace{ \matr{\Gamma}_{ m q_1 } \dots \matr{\Gamma}_{ q_{ l - 1 } n } }_{\mbox{$l$ terms}}$$
and $\tr ( \matr{\Gamma} ^ l )$ by
\begin{eqnarray*}
 \tr ( \matr{\Gamma} ^ l ) & = & \sum_{ n = 1 } ^ N \tr \left[ ( \matr{\Gamma} ^ l )_{ n n } \right] = \sum_{ n = 1 } ^ N \tr \left( \sum_{ q_1, \dots, q_{ l - 1 } } \matr{\Gamma}_{ n q_{ l - 1 } } \dots \matr{\Gamma}_{ q_1 n } \right).
\end{eqnarray*}
Consider then the directed and fully connected graph with $N$ nodes $\{ 1, \dots, N \}$, an arrow from any $m$ to any $n \neq m$ (no self-connections), and corresponding (potentially matrix) weight $\matr{\Gamma}_{ n | m }$. In this graph, a directed loop is a directed path that begins and starts at the same node. It is a $k$-loop if the directed path is composed of exactly $k$ arrows. For any node $n$ and integer $k$, let $\mathcal{L}_k ^ \rightarrow ( n )$ be the set of all directed $k$-loops starting and ending at node $n$ and $\mathcal{L}_k ^ \rightarrow$ the set of all $k$-loops. For any directed path $p = ( q_1 \to \dots \to q_k )$, define $\tau ( p )$ as the trace of the product of the weights along $p$
$$\tau ( p ) = \tr \left( \matr{\Gamma}_{ q_k | q_{ k - 1 } } \dots \matr{\Gamma}_{ q_2 | q_1 } \right).$$
With these notations, $\tr ( \matr{\Gamma} ^ l )$ can be interpreted as the sum of the values taken by $\tau$ along all directed $l$-loops starting and ending at every node of the graph
$$\tr \left( \matr{\Gamma} ^ l \right) = \sum_{ n = 1 } ^ N \sum_{ p \in \mathcal{L}_l ^ \rightarrow ( n ) } \tau ( p ).$$
It can also be seen as the sum of the values taken by $\tau$ along all directed $l$-loops
$$\tr \left( \matr{\Gamma} ^ l \right) = \sum_{ p \in \mathcal{L}_l ^ \rightarrow } \tau ( p ).$$
See Figure~\ref{fig} for an illustration of this interpretation. Note that the fact that 1-loops do not exist is interpreted as $\tr ( \matr{\Gamma} ) = 0$.
\par
This graphical interpretation is in particlar compatible with a partitioning into two subvectors. In that case, the corresponding directed graph only has two nodes and two arrows, and $\mathcal{L}_l ^ \rightarrow = \emptyset$ for $l$ odd, which is in agreement with the fact that all cumulants of odd order are equal to zero.

\section{Discussion} \label{s:disc}

In the present manuscript, we introduced multiinformation density, a random variable that generalizes information density and whose expectation defines multiinformation. We focused on the case of a multivariate normal distribution and derived a closed form for its cumulant-generating function as well as its cumulants. We showed that the cumulant-generating function does not depend on the values taken by the variance coefficients of the covariance matrix, and that the computations required have a simple graphical interpretation. We also considered the special case of a partitioning into two subvectors, showing the relationship between our results and existing quantities.
\par
Interestingly, the results show that the cumulant-generating function of multiinformation density is a function of a specific quantity, namely the block matrix $\matr{\Gamma}$ composed of all matrices of regression coefficients corresponding to the partitioning of $\vect{X}$ into $( \vect{X}_1, \dots, \vect{X}_N )$ as defined in Equation~\eqref{eq:def:regcoeff}. This entails that the probability distribution of multiinformation density as well as all its moments are fully defined by $\matr{\Gamma}$. In particular, it can be shown that mutual information and multiinformation themselves, as expectations, are functions of $\matr{\Gamma}$ only, namely (see Appendix~\ref{an:im})
$$I ( \vect{X}_1; \dots; \vect{X}_N ) = - \frac{ 1 } { 2 } \ln | \matr{I}_d + \matr{\Gamma} |.$$
Going back to our question of knowing whether the multivariate normal distribution family could be parameterized in a natural way with emphasis on its dependence properties, it can be argued that $\matr{\Gamma}$ is a good candidate to this aim.

\paragraph{Generalization}

The expression of the cumulant-generating function of multiinformation in the case of multivariate normal distributions is quite simple. This simplicity is mostly a consequence of the stability of the multivariate normal family to most operations performed here: (i) the product of the marginals is also a multivariate normal distribution (this is strongly related to the fact that independence and uncorrelatedness are equivalent for multivariate normal distributions); (ii) the ratio of the joint distribution to its marginals takes a simple form that is again closely related to the multivariate normal family; and (iii) the exponentiation of $t$ times the log of this ratio still has the form of a multivariate normal distribution. It would be of interest to determine such a simplicity would still hold in more general settings such as more general distribution families and more general functions of the ratio $f_{ \vect{X} } ( \vect{X} ) / \prod f_{ \vect{X}_i } ( \vect{X}_i )$.
\par
Regarding the type of family considered, a first step could be to consider the family of multivariate $t$ distributions \gcite{Kotz-2004}. In this case, $f_{ \vect{X} } ( \vect{x} )$ and $f_{ \vect{Y} } ( \vect{y} )$ are both multivariate $t$ distributions \gcite[\S1.10]{Kotz-2004} but $f_{ \vect{X} } ( \vect{x} ) \, f_{ \vect{Y} } ( \vect{y} )$ itself is not a multivariate $t$ distribution, and the cumulant-generating function of the log of the resulting ratio does not have a simple closed form. Another, further step would be to consider elliptically contoured distributions \gcite[\S2.7]{Anderson_TW-2003}, but again no simplification seems to occur.
\par
Another potential generalization would be obtained by replacing the log function of multiinformation density with a quantity of the form
$$g \left[ \frac{ f_{ \vect{X} } ( \vect{X} ) } { \prod f_{ \vect{X}_i } ( \vect{X}_i ) } \right],$$
in the spirit of $f$-divergences \gcite{Csiszar-1963, Csiszar-1967, Ali-1966, Vajda-1972}. However, by replacing $g ( t ) = \ln ( t )$ with another, more general form, we lose the key property that the log is the inverse of the exponentiation taken to compute the cumulant-generating function, which also simplifies calculations.

\paragraph{Future work}

Investigating the level of dependence between variables is still a thorny issue. Such an issue is usually tackled by investigating the properties of mutual information or multiinformation. We believe that considering multiinformation density instead of multiinformation (i.e., the random variable instead of its expectation) couldt contribute to a better characterization of dependence. For instance, the variance of $i_d$, given in Equation~(\ref{eq:idep:var}), could be considered, in addition to or instead of multiinformation, to quantify the presence or absence of dependence. An important point to advance in this direction would be to provide estimators for the quantities obtained here. In the case of two variables, an estimator of this quantity, expressed in Equation~(\ref{eq:2:var}), is readily available \gcite{Jupp-1980}. Its generalization to more than two variables could yield a new tool to investigate mutual independence.
\par
Another point of interest is the investigation of the behavior of multiinformation density in the case of large dimension. For instance, the asymptotic normality of this measure could be of interest. A quick derivation (see Appendix~\ref{an:norm}) shows that, in the simple case of a homogeneous correlation matrix and a partitioning into 1-dimensional vectors, $i_d$ is not asymptotically normal. We hope that tackling this issue in the more general case will help better understand dependences within systems composed of many variables.

\appendix

\section{Mutual independence, $i_d$ and multiinformation} \label{ann:id}

We here show the equivalence between the following
\begin{enumerate}
 \item The $\vect{X}_n$'s are mutually independent;
 \item $i_d \equiv 0$;
 \item $I ( \vect{X}_1; \dots; \vect{X}_N ) = 0$;
 \item $\var ( i_d ) = 0$.
\end{enumerate}
The equivalence between (1) and (2) is straightforward by definition of $i_d$. If $i_d \equiv 0$, then all its moments of $i_d$ are equal to 0, including its expectation (multiinformation) and its variance, leading to $I ( \vect{X}_1; \dots; \vect{X}_N ) = 0$ and $\var ( i_d ) = 0$.
\par
Since multiinformation is a Kullback-Leibler divergence, $I ( \vect{X}_1; \dots; \vect{X}_N ) = 0$, entails that we have
$$f ( \vect{x} ) = \prod_{ n = 1 } ^ N f_n ( \vect{x}_n )$$
for all $\vect{x}$, i.e., the  $\vect{X}_n$'s are mutually independent.
\par
Finally, if $\var ( i_d ) = 0$, then $i_d$ is a constant, that is,
$$f ( \vect{X} ) = k \prod_{ n = 1 } ^ N f_n ( \vect{x}_n ).$$
The fact that $f$ and the $f_n$'s are distributions, and therefore must norm to 1, entails that we must have $k = 1$, that is, $i_d \equiv 0$.

\section{Positive-definiteness of $\matr{\Sigma} ^ { - 1 } - t \matr{\Phi}$} \label{an:pd}

We need to show that $\matr{\Sigma} ^ { - 1 } - t \matr{\Phi}$ is a symmetric positive definite matrix in a neighborhood of $t = 0$. This matrix can be expressed as
$$\matr{\Sigma} ^ { - 1 } - t \matr{\Phi} = ( 1 + t ) \matr{\Sigma} ^ { - 1 } - t \, \mathrm{diag} ( \matr{\Sigma}_{ 1 1 }, \dots, \matr{\Sigma}_{ N N } ) ^ { - 1 }.$$
As a difference of two symmetric matrices, it is also symmetric. Furthermore, since the two matrices in the right-hand side of the equation are positive definite they are diagonalizable in the same basis, i.e., there exists a nonsingular matrix $\matr{F}$ such that \gcite[Theorem~A.2.2]{Anderson_TW-2003}
$$\matr{F} \transp \matr{\Sigma} ^ { - 1 } \matr{F} = \begin{pmatrix}
 \lambda_1^2 & & \\
 & \ddots & \\
 & & \lambda_d^2
\end{pmatrix}$$
and
$$\matr{F} \transp \mathrm{diag} ( \matr{\Sigma}_{ 1 1 }, \dots, \matr{\Sigma}_{ N N } ) ^ { - 1 } \matr{F} = \matr{I}.$$
Since $\matr{\Sigma} ^ { - 1 }$ is positive definite, we furthermore have $\lambda_i ^ 2 > 0$. $\matr{\Sigma} ^ { - 1 } - t \matr{\Phi}$ is therefore diagonalizable as well, with eigenvalues given by $(1+t) \lambda_i^2-t = (\lambda_i^2-1)t+\lambda_i^2$, which is positive in a neighborhood of $t = 0$. $\matr{\Sigma} ^ { - 1 } - t \matr{\Phi}$ is therefore positive definite in a neighborhood of $t = 0$.


\section{Alternative expression of multiinformation} \label{an:im}

For a decomposition of a multidimensional normal variable into several subvectors, multiinformation reads
$$I ( \vect{X}_1; \dots; \vect{X}_N ) = \frac{ 1 } { 2 } \ln \frac{ \prod_{ n = 1 } ^ N | \matr{\Sigma}_{ n n } | } { | \matr{\Sigma} | }.$$
By comparison, we calculate
\begin{eqnarray*}
 \matr{I}_d + \matr{\Gamma} & = & \matr{I}_d + \matr{\Sigma} \matr{\Phi} \\
 & = & \matr{I}_d + \matr{\Sigma} \left[  \mathrm{diag} ( \matr{\Sigma}_{ 1 1 }, \dots, \matr{\Sigma}_{ N N } ) ^ { - 1 } - \matr{\Sigma} ^ { - 1 } \right] \\
 & = & \matr{\Sigma} \, \mathrm{diag} ( \matr{\Sigma}_{ 1 1 }, \dots, \matr{\Sigma}_{ N N } ) ^ { - 1 },
\end{eqnarray*}
leading to
\begin{eqnarray*}
 | \matr{I}_d + \matr{\Gamma} | & = & | \matr{\Sigma} \, \mathrm{diag} ( \matr{\Sigma}_{ 1 1 }, \dots, \matr{\Sigma}_{ N N } ) ^ { - 1 } | \\
 & = & \frac { | \matr{\Sigma} | } { \prod_{ n = 1 } ^ N | \matr{\Sigma}_{ n n } | },
\end{eqnarray*}
and, finally,
$$- \frac{1}{2} \ln | \matr{I}_d + \matr{\Gamma} | = \frac{1}{2} \ln \frac { \prod_{ n = 1 } ^ N | \matr{\Sigma}_{ n n } | } { | \matr{\Sigma} | }.$$

\section{Checking asymptotic normality} \label{an:norm}

Let the correlation matrix $\matr{R}_d$ be a $d$-by-$d$ homogeneous matrix with parameter $\rho$, i.e., a matrix with 1s on the diagonal and all off-diagonal elements equal to $\rho$. such a matrix has two eigenvalues: $1 + ( d - 1 ) \rho$ with multiplicity 1 (associated with the vector composed only of 1s) and $1 - \rho$ with multiplicity $d - 1$ (associated with the subspace of vectors with a zero mean).  Such a matrix is positive definite for
$$ - \frac{1}{d-1} \leq \rho < 1.$$
The expectation of $i_d$ is given by
$$- \frac{ 1 } { 2 } \left\{ ( d - 1 ) \ln ( 1 - \rho ) + \ln [ 1 + ( d - 1 ) \rho ] \right\}$$
To compute the higher cumulants of $i_d$, let $\matr{U}_d$ the $d$-by-$d$ matrix with all elements equal to 1. Using the fact that $\matr{\Gamma} = \rho ( \matr{U}_d - \matr{I}_d )$ together with $\matr{U}_d ^ l = d ^ { l - 1 } \matr{U}_d$ for $l \geq 2$, we obtain
\begin{eqnarray*}
 \matr{\Gamma} ^ l & = & \rho ^ l \left[ ( - 1 ) ^ l \matr{I}_d + \frac{ ( d - 1 ) ^ l - ( - 1 ) ^ l } { d } \matr{U}_d \right] \\
 \tr ( \matr{\Gamma} ^ l ) & = & \rho ^ l \left[ ( - 1 ) ^ l d + ( d - 1 ) ^ l - ( - 1 ) ^ l \right] \\
 \kappa_l ( i_d ) & = & \frac{ ( l - 1 ) ! } { 2 }  \rho ^ l \left[ ( - 1 ) ^ l d + ( d - 1 ) ^ l - ( - 1 ) ^ l \right].                                                                                                                                                                                                              \end{eqnarray*}
In particular, we have $\var( i_d ) = \rho ^ 2 d ( d - 1 ) / 2$. For large $d$, we have $\kappa_l ( i_d ) \sim \frac{ ( l - 1 ) ! } { 2 }  \rho ^ l d ^ l$ for $l \geq 2$ and, in particular, $\var( i_d ) \sim \rho ^ 2 d ^ 2 / 2$. To investigate the asymptotic normality of $i_d$, we classically consider $u = [ i_d - \esp ( i_d ) ] / \sqrt{ \var ( i_d ) }$. Using the fact that the cumulant of order $l$ is homogeneous of degree $l$, we obtain $\kappa_l ( u ) \sim 2 ^ { l/2 - 1 } ( l -1 ) ! = \mbox{cste}$. If $u$ were asymptotically normal, $\kappa_l ( u )$ for $l \geq 3$ would tend to 0 as $d \to \infty$, which is not the case. As a consequence, $u$ is not asymptotically normal.

% \bibliographystyle{elsarticle-harv}
% \bibliography{nonabrev,anglais,mabiblio}

\begin{thebibliography}{18}
\expandafter\ifx\csname natexlab\endcsname\relax\def\natexlab#1{#1}\fi
\providecommand{\url}[1]{\texttt{#1}}
\providecommand{\href}[2]{#2}
\providecommand{\path}[1]{#1}
\providecommand{\DOIprefix}{doi:}
\providecommand{\ArXivprefix}{arXiv:}
\providecommand{\URLprefix}{URL: }
\providecommand{\Pubmedprefix}{pmid:}
\providecommand{\doi}[1]{\href{http://dx.doi.org/#1}{\path{#1}}}
\providecommand{\Pubmed}[1]{\href{pmid:#1}{\path{#1}}}
\providecommand{\bibinfo}[2]{#2}
\ifx\xfnm\relax \def\xfnm[#1]{\unskip,\space#1}\fi
%Type = Book
\bibitem[{Abramowitz and Stegun(1972)}]{Abramowitz-1972}
\bibinfo{editor}{Abramowitz, M.}, \bibinfo{editor}{Stegun, I.A.} (Eds.),
  \bibinfo{year}{1972}.
\newblock \bibinfo{title}{Handbook of Mathematical Functions}.
\newblock Number~\bibinfo{number}{55} in \bibinfo{series}{Applied Math.},
  \bibinfo{publisher}{National Bureau of Standards}.
%Type = Article
\bibitem[{Ali and Silvey(1966)}]{Ali-1966}
\bibinfo{author}{Ali, S.}, \bibinfo{author}{Silvey, S.D.},
  \bibinfo{year}{1966}.
\newblock \bibinfo{title}{A general class of coefficients of divergence of one
  distribution from another}.
\newblock \bibinfo{journal}{Journal of the Royal Statistical Society: Series B
  (Statistical Methodology)} \bibinfo{volume}{28}, \bibinfo{pages}{131--142}.
%Type = Book
\bibitem[{Anderson(2003)}]{Anderson_TW-2003}
\bibinfo{author}{Anderson, T.W.}, \bibinfo{year}{2003}.
\newblock \bibinfo{title}{An Introduction to Multivariate Statistical
  Analysis}.
\newblock Wiley Series in Probability and Mathematical Statistics.
  \bibinfo{edition}{3rd} ed., \bibinfo{publisher}{John Wiley and Sons, New
  York}.
%Type = Article
\bibitem[{Csisz{\'a}r(1963)}]{Csiszar-1963}
\bibinfo{author}{Csisz{\'a}r, I.}, \bibinfo{year}{1963}.
\newblock \bibinfo{title}{Eine informationstheoretische {U}ngleichung und ihre
  {A}nwendung auf den {B}eweis der {E}rgodizitat von {M}arkoffschen {K}etten}.
\newblock \bibinfo{journal}{A Magyar Tudom{\'a}nyos Akad{\'e}mia Matematikai
  {\'e}s Fizikai Tudom{\'a}nyok Oszt{\'a}ly{\'a}nak K{\"o}zlem{\'e}nyei}
  \bibinfo{volume}{8}, \bibinfo{pages}{85--108}.
%Type = Article
\bibitem[{Csisz{\'a}r(1967)}]{Csiszar-1967}
\bibinfo{author}{Csisz{\'a}r, I.}, \bibinfo{year}{1967}.
\newblock \bibinfo{title}{Information-type measures of difference of
  probability distributions and indirect observation}.
\newblock \bibinfo{journal}{Studia Scientiarum Mathematicarum Hungarica}
  \bibinfo{volume}{2}, \bibinfo{pages}{229--318}.
%Type = Book
\bibitem[{Garner(1962)}]{Garner-1962}
\bibinfo{author}{Garner, W.R.}, \bibinfo{year}{1962}.
\newblock \bibinfo{title}{Uncertainty and Structure as Psychological Concepts}.
\newblock \bibinfo{publisher}{John Wiley \& Sons, New York}.
%Type = Incollection
\bibitem[{Higham(2007)}]{Higham-2007}
\bibinfo{author}{Higham, N.J.}, \bibinfo{year}{2007}.
\newblock \bibinfo{title}{Functions of matrices}, in: \bibinfo{editor}{Hogben,
  L.} (Ed.), \bibinfo{booktitle}{Handbook of Linear Algebra}.
  \bibinfo{publisher}{Chapman \& Hall/CRC Press, Boca Raton}. Discrete
  Mathematics and its Applications. chapter~\bibinfo{chapter}{11}.
%Type = Article
\bibitem[{Joe(1989)}]{Joe-1989b}
\bibinfo{author}{Joe, H.}, \bibinfo{year}{1989}.
\newblock \bibinfo{title}{Relative entropy measures of multivariate
  dependence}.
\newblock \bibinfo{journal}{Journal of the American Statistical Association}
  \bibinfo{volume}{84}, \bibinfo{pages}{157--164}.
%Type = Article
\bibitem[{Jupp and Mardia(1980)}]{Jupp-1980}
\bibinfo{author}{Jupp, P.E.}, \bibinfo{author}{Mardia, K.V.},
  \bibinfo{year}{1980}.
\newblock \bibinfo{title}{A general correlation coefficient for directional
  data and related regression problems}.
\newblock \bibinfo{journal}{Biometrika} \bibinfo{volume}{67},
  \bibinfo{pages}{163--173}.
%Type = Book
\bibitem[{Kendall(1945)}]{Kendall-1945}
\bibinfo{author}{Kendall, M.G.}, \bibinfo{year}{1945}.
\newblock \bibinfo{title}{The Advanced Theory of Statistics}.
  volume~\bibinfo{volume}{1}.
\newblock \bibinfo{edition}{2nd} ed., \bibinfo{publisher}{Charles Griffin \&
  Co. Ltd., London}.
%Type = Book
\bibitem[{Kotz and Nadarajah(2004)}]{Kotz-2004}
\bibinfo{author}{Kotz, S.}, \bibinfo{author}{Nadarajah, S.},
  \bibinfo{year}{2004}.
\newblock \bibinfo{title}{Multivariate $t$ Distributions and their
  Applications}.
\newblock \bibinfo{publisher}{Cambridge University Press, Cambridge, UK}.
%Type = Book
\bibitem[{Kullback(1968)}]{Kullback-1968}
\bibinfo{author}{Kullback, S.}, \bibinfo{year}{1968}.
\newblock \bibinfo{title}{Information Theory and Statistics}.
\newblock \bibinfo{publisher}{Dover, Mineola, NY}.
%Type = Article
\bibitem[{Magnus(1986)}]{Magnus-1986}
\bibinfo{author}{Magnus, J.}, \bibinfo{year}{1986}.
\newblock \bibinfo{title}{The exact moments of a ratio of quadratic forms}.
\newblock \bibinfo{journal}{Annales d'{\'e}conomie et de statistique}
  \bibinfo{volume}{4}, \bibinfo{pages}{95--109}.
%Type = Book
\bibitem[{Mardia and Jupp(2000)}]{Mardia-2000}
\bibinfo{author}{Mardia, K.V.}, \bibinfo{author}{Jupp, P.E.},
  \bibinfo{year}{2000}.
\newblock \bibinfo{title}{Directional Statistics}.
\newblock Wiley Series in Probability and Statistics,
  \bibinfo{publisher}{Wiley, Chichester}.
%Type = Misc
\bibitem[{Polyanskiy and Wu(2017)}]{Polyanskiy-2017}
\bibinfo{author}{Polyanskiy, Y.}, \bibinfo{author}{Wu, Y.},
  \bibinfo{year}{2017}.
\newblock \bibinfo{title}{Lecture notes on information theory}.
\newblock \bibinfo{howpublished}{http://www.stat.yale.edu/$\sim$yw562/ln.html}.
%Type = Inproceedings
\bibitem[{Studen{\'y} and Vejnarov{\'a}(1998)}]{Studeny-1998b}
\bibinfo{author}{Studen{\'y}, M.}, \bibinfo{author}{Vejnarov{\'a}, J.},
  \bibinfo{year}{1998}.
\newblock \bibinfo{title}{The multiinformation function as a tool for measuring
  stochastic dependence}, in: \bibinfo{editor}{Jordan, M.I.} (Ed.),
  \bibinfo{booktitle}{Proceedings of the NATO Advanced Study Institute on
  Learning in Graphical Models}, pp. \bibinfo{pages}{261--298}.
%Type = Article
\bibitem[{Vajda(1972)}]{Vajda-1972}
\bibinfo{author}{Vajda, I.}, \bibinfo{year}{1972}.
\newblock \bibinfo{title}{On the $f$-divergence and singularity of probability
  measures}.
\newblock \bibinfo{journal}{Periodica Mathematica Hungarica}
  \bibinfo{volume}{2}, \bibinfo{pages}{223--234}.
%Type = Article
\bibitem[{Watanabe(1960)}]{Watanabe-1960}
\bibinfo{author}{Watanabe, S.}, \bibinfo{year}{1960}.
\newblock \bibinfo{title}{Information theoretical analysis of multivariate
  correlation}.
\newblock \bibinfo{journal}{IBM Journal of Research and Development}
  \bibinfo{volume}{4}, \bibinfo{pages}{66--82}.

\end{thebibliography}

\end{document}